\documentclass[11pt,a4paper]{article}

\usepackage{amsmath,amssymb,amsthm}

\usepackage[T1]{fontenc}
\usepackage[applemac]{inputenc}	
\usepackage[english]{babel}
\usepackage{geometry}

\newtheorem{theorem}{Theorem}
\newtheorem{corollary}[theorem]{Corollary}
\newtheorem{lemma}[theorem]{Lemma}

\newtheorem{conjecture}[theorem]{Conjecture}
 
 \theoremstyle{remark}

\newcommand{\comment}[1]{}

\newcommand{\R}{\mathbb R}

\newcommand{\eps}{\varepsilon}

\newcommand{\de}{{\rm d}}
\newcommand{\supp}{{\rm supp}}

 \newgeometry{tmargin=2cm, bmargin=2cm, lmargin=2cm, rmargin=2cm}

\begin{document}
\large \title{Borell's generalized Pr\'ekopa-Leindler inequality: A simple proof}

\author{Arnaud Marsiglietti\thanks{Supported in part by the Institute for Mathematics and its Applications with funds provided by the National Science Foundation.}}

\date{}

\maketitle
\begin{abstract}

We present a simple proof of Christer Borell's general inequality in the Brunn-Minkowski theory. We then discuss applications of Borell's inequality to the log-Brunn-Minkowski inequality of B\"or\"oczky, Lutwak, Yang and Zhang.

\end{abstract}

\noindent {\bf 2010 Mathematics Subject Classification.}  Primary 28A75, 52A40.

\noindent {\bf Keywords.} Brunn-Minkowski, Convex body, log-Brunn-Minkowski inequality, mass transportation.

\section{Introduction}

Let us denote by $\supp(f)$ the support of a function $f$. In~\cite{B2} Christer Borell proved the following inequality (see \cite[Theorem 2.1]{B2}), which we will call the Borell-Brunn-Minkowski inequality.

\begin{theorem}[Borell-Brunn-Minkowski inequality]\label{Borell}

Let $f,g,h : \R^n \to [0, +\infty)$ be measurable functions. Let $\varphi = (\varphi_1, \dots, \varphi_n) : \supp(f) \times \supp(g) \to \R^n$ be a continuously differentiable function with positive partial derivatives, such that $\varphi_k(x,y) = \varphi_k(x_k,y_k)$ for every $x=(x_1, \dots, x_n) \in \supp(f)$, $y=(y_1, \dots, y_n) \in \supp(g)$. Let $\Phi : [0, +\infty) \times [0, +\infty) \to [0, + \infty)$ be a continuous function, homogeneous of degree 1 and increasing in each variable. If the inequality
\begin{eqnarray}\label{condition}
h(\varphi(x,y)) \Pi_{k=1}^n \left( \frac{\partial \varphi_k}{\partial x_k} \rho_k + \frac{\partial \varphi_k}{\partial y_k} \eta_k \right) \geq \Phi(f(x)\Pi_{k=1}^n \rho_k, g(y) \Pi_{k=1}^n \eta_k)
\end{eqnarray}
holds for every $x \in \supp(f)$, for every $y \in \supp(g)$, for every $\rho_1, \dots, \rho_n > 0$ and for every $\eta_1, \dots, \eta_n > 0$, then
$$ \int h \geq \Phi \left(\int f, \int g \right). $$

\end{theorem}

C. Borell proved a slightly more general statement, involving an arbitrary number of functions. For simplicity, we restrict ourselves to the statement of Theorem~\ref{Borell}.

Theorem \ref{Borell} yields several important consequences. For example, applying Theorem \ref{Borell} to indicators of compact sets (i.e. $f=1_A$, $g=1_B$, $h=1_{\varphi(A,B)}$) yields the following generalized Brunn-Minkowski inequality.

\begin{corollary}[Generalized Brunn-Minkowski inequality]\label{GBM}

Let $A,B$ be compact subsets of $\R^n$. Let $\varphi = (\varphi_1, \dots, \varphi_n) : A \times B \to \R^n$ be a continuously differentiable function with positive partial derivatives, such that $\varphi_k(x,y) = \varphi_k(x_k,y_k)$ for every $x=(x_1, \dots, x_n) \in A$, $y=(y_1, \dots, y_n) \in B$. Let $\Phi : [0, +\infty) \times [0, +\infty) \to [0, + \infty)$ be a continuous function, homogeneous of degree~1 and increasing in each variable. If the inequality
$$ \Pi_{k=1}^n \left( \frac{\partial \varphi_k}{\partial x_k} \rho_k + \frac{\partial \varphi_k}{\partial y_k} \eta_k \right) \geq \Phi(\Pi_{k=1}^n \rho_k, \Pi_{k=1}^n \eta_k) $$
holds for every $\rho_1, \dots, \rho_n, \eta_1, \dots, \eta_n > 0$, then
$$ |\varphi(A,B)| \geq \Phi \left(|A|, |B| \right), $$
where $|\cdot|$ denotes Lebesgue measure and $\varphi(A,B) = \{ \varphi(x,y) : x \in A, y \in B \}$.

\end{corollary}

The classical Brunn-Minkowski inequality (see e.g. \cite{Sch}, \cite{G}) follows from Corollary \ref{GBM} by taking $\varphi(x,y) = x+y$, $x \in A, y \in B$, and $\Phi(a,b) = (a^{1/n} + b^{1/n})^n$, $a,b \geq 0$. Although the Brunn-Minkowski inequality goes back to more than a century ago, it still attracts a lot of attention (see e.g. \cite{NT}, \cite{CN}, \cite{GHW}, \cite{M3}, \cite{CM}, \cite{CDP}, \cite{CGN}, \cite{GS}, \cite{M5}).

Theorem \ref{Borell} also allows us to recover the so-called Borell-Brascamp-Lieb inequality. Let us denote by $M_s^{\lambda}(a,b)$ the $s$-mean of the real numbers $a, b \geq 0$ with weight $\lambda \in [0,1]$, defined as
$$ M_s^{\lambda}(a,b) = ((1-\lambda)a^s + \lambda b^s)^{\frac{1}{s}} \quad \mbox{if $s \notin \{-\infty, 0, +\infty\}$}, $$
$M_{-\infty}^{\lambda}(a,b) = \min(a,b)$, $M_0^{\lambda}(a,b) = a^{1-\lambda} b^{\lambda}$, $M_{+\infty}^{\lambda}(a,b) = \max(a,b)$. We will need the following H\"older inequality (see e.g. \cite{HLP}).

\begin{lemma}[Generalized H\"older inequality]\label{holder}

Let  $\alpha, \beta, \gamma \in \R \cup \{ +\infty \}$ such that $\beta + \gamma \geq 0$ and $\frac{1}{\beta} + \frac{1}{\gamma} = \frac{1}{\alpha}$. Then, for every $a,b,c,d \geq 0$ and $\lambda \in [0,1]$,
\begin{eqnarray*}
M_{\alpha}^{\lambda}(ac, bd) \leq M_{\beta}^{\lambda}(a,b) M_{\gamma}^{\lambda}(c,d).
\end{eqnarray*}

\end{lemma}

\begin{corollary}[Borell-Brascamp-Lieb inequality]\label{BBL}

Let $\gamma \geq -\frac{1}{n}$, $\lambda \in [0,1]$ and $f,g,h : \R^n \to [0, +\infty)$ be measurable functions. If the inequality
$$ h((1-\lambda)x + \lambda y) \geq M_{\gamma}^{\lambda}(f(x), g(y)) $$
holds for every $x \in \supp(f), y \in \supp(g)$, then
$$ \int_{\R^n} h \geq M_{\frac{\gamma}{1+\gamma n}}^{\lambda}\left(\int_{\R^n} f, \int_{\R^n} g\right). $$

\end{corollary}

Corollary \ref{BBL} follows from Theorem \ref{Borell} by taking $\varphi(x,y) = (1-\lambda)x + \lambda y$, $x \in \supp(f), y \in \supp(g)$, and $\Phi(a,b) = M_{\frac{\gamma}{1+\gamma n}}^{\lambda}(a, b)$, $a,b \geq 0$. Indeed, using Lemma \ref{holder}, one obtains that for every $x \in \supp(f), y \in \supp(g)$, and for every $\rho_1, \dots, \rho_n, \eta_1, \dots, \eta_n > 0$,
\begin{eqnarray*}
h(\varphi(x,y)) \Pi_{k=1}^n \left( \frac{\partial \varphi}{\partial x_k} \rho_k + \frac{\partial \varphi}{\partial y_k} \eta_k \right) & = & h((1-\lambda)x + \lambda y) \Pi_{k=1}^n ((1-\lambda) \rho_k + \lambda \eta_k) \\ & \geq & M_{\gamma}^{\lambda}(f(x), g(y)) M_{\frac{1}{n}}^{\lambda}(\Pi_{k=1}^n \rho_k, \Pi_{k=1}^n \eta_k) \\ & \geq &  M_{\frac{\gamma}{1+\gamma n}}^{\lambda}(f(x)\Pi_{k=1}^n \rho_k, g(y)\Pi_{k=1}^n \eta_k) \\ & = & \Phi(f(x)\Pi_{k=1}^n \rho_k, g(y) \Pi_{k=1}^n \eta_k).
\end{eqnarray*}

Corollary \ref{BBL} was independently proved by Borell (see \cite[Theorem 3.1]{B2}), and by Brascamp and Lieb \cite{BL}.

Another important consequence of the Borell-Brunn-Minkowski inequality is obtained when considering $\varphi$ to be nonlinear. Let us denote for ${\bf p}=(p_1, \dots, p_n) \in [- \infty, +\infty]^n$, $x=(x_1, \dots, x_n) \in [0, +\infty]^n$ and $y=(y_1, \dots, y_n) \in [0, +\infty]^n$,
$$ M_{\bf p}^{\lambda}(x,y) = (M_{p_1}^{\lambda}(x_1, y_1), \dots, M_{p_n}^{\lambda}(x_n, y_n)). $$

\begin{corollary}[nonlinear extension of the Brunn-Minkowski inequality]\label{B2}

Let ${\bf p}=(p_1, \dots, p_n) \in [0,1]^n$, $\gamma \geq -(\sum_{i=1}^n p_i^{-1})^{-1}$, $\lambda \in [0,1]$, and $f,g,h : [0,+\infty)^n \to [0, +\infty)$ be measurable functions. If the inequality
$$ h(M_{\bf p}^{\lambda}(x,y)) \geq M_{\gamma}^{\lambda}(f(x),g(y)) $$
holds for every $x \in \supp(f), y \in \supp(g)$, then
$$ \int_{[0,+\infty)^n} h \geq M_{(\sum_{i=1}^n{p_i^{-1}} + \gamma^{-1})^{-1}}^{\lambda}\left( \int_{[0,+\infty)^n} f, \int_{[0,+\infty)^n} g \right). $$

\end{corollary}

Corollary \ref{B2} follows from Theorem \ref{Borell} by taking $\varphi(x,y) = M_{\bf p}^{\lambda}(x,y)$, $x \in \supp(f), y \in \supp(g)$, and $\Phi(a,b) = M_{(\sum_{i=1}^n{p_i^{-1}} + \gamma^{-1})^{-1}}^{\lambda}(a, b)$, $a,b \geq 0$. Indeed, using Lemma \ref{holder}, one obtains that for every $x \in \supp(f), y \in \supp(g)$, and for every $\rho_1, \dots, \rho_n, \eta_1, \dots, \eta_n > 0$,
\begin{eqnarray*}
h(\varphi(x,y)) \Pi_{k=1}^n \left( \frac{\partial \varphi}{\partial x_k} \rho_k + \frac{\partial \varphi}{\partial y_k} \eta_k \right) & = & h(M_{\bf p}^{\lambda}(x,y)) \Pi_{k=1}^n M_{\frac{p_k}{1-p_k}}^{\lambda}(x_k^{1-p_k}, y_k^{1-p_k}) M_1(x_k^{p_k-1}\rho_k, y_k^{p_k-1} \eta_k) \\ & \geq & M_{\gamma}^{\lambda}(f(x), g(y)) \Pi_{k=1}^n M_{p_k}^{\lambda}(\rho_k, \eta_k) \\ & \geq &  M_{\gamma}^{\lambda}(f(x), g(y)) M_{(\sum_{i=1}^n p_i^{-1})^{-1}}^{\lambda}(\Pi_{k=1}^n \rho_k, \Pi_{k=1}^n \eta_k) \\ & \geq & M_{(\sum_{i=1}^n{p_i^{-1}} + \gamma^{-1})^{-1}}^{\lambda}(f(x)\Pi_{k=1}^n \rho_k, g(y) \Pi_{k=1}^n \eta_k) \\ & = & \Phi(f(x)\Pi_{k=1}^n \rho_k, g(y) \Pi_{k=1}^n \eta_k).
\end{eqnarray*}

In the particular case where ${\bf p}=(0, \dots, 0)$, Corollary \ref{B2} was rediscovered by Ball \cite{Ball}. In the general case, Corollary \ref{B2} was rediscovered by Uhrin \cite{U}. 

Notice that the condition on $p$ in Corollary \ref{B2} is less restrictive in dimension 1. It reads as follows:

\begin{corollary}[nonlinear extension of the Brunn-Minkowski inequality on the line]\label{line}

Let $p \leq 1$, $\gamma \geq -p$, and $\lambda \in [0,1]$. Let $f,g,h : [0,+\infty) \to [0, +\infty)$ be measurable functions such that for every $x \in \supp(f), y \in \supp(g)$,
$$ h(M_p^{\lambda}(x,y)) \geq M_{\gamma}^{\lambda}(f(x),g(y)). $$
Then,
$$ \int_0^{+\infty} h \geq M_{\left(\frac{1}{p} + \frac{1}{\gamma}\right)^{-1}}^{\lambda}\left( \int_0^{+\infty} f, \int_0^{+\infty} g \right). $$

\end{corollary}

A simple proof of Corollary \ref{line} was recently given by Bobkov et al. \cite{BCF}. \\

In section 2, we present a simple proof of Theorem \ref{Borell}, based on mass transportation. In section 3, we discuss applications of the above inequalities to the log-Brunn-Minkowski inequality of B\"or\"oczky, Lutwak, Yang and Zhang. We also prove an equivalence between the log-Brunn-Minkowski inequality and its possible extensions to convex measures (see section 3 for precise definitions).

\section{A simple proof of the Borell-Brunn-Minkowski inequality}

In this section, we present a simple proof of Theorem \ref{Borell}.

\begin{proof}[Proof of Theorem \ref{Borell}]
The proof is done by induction on the dimension. To prove the theorem in dimension 1, we use a mass transportation argument. \\

{\bf Step 1 :} (In dimension 1) \\
First let us see that if $\int f = 0$ or $\int g = 0$, then the result holds. Let us assume, without loss of generality, that $\int g = 0$. By taking $\rho=1$, by letting $\eta$ go to $0$ and by using continuity and homogeneity of $\Phi$ in the condition (\ref{condition}), one obtains
$$ h(\varphi(x,y)) \frac{\partial \varphi}{\partial x} \geq \Phi(f(x),0) = f(x)\Phi(1,0). $$
It follows that, for fixed $y \in \supp(g)$,
$$ \int h(z) \de z \geq \int_{\varphi(\supp(f),y)} h(z) \de z = \int_{\supp(f)} h(\varphi(x,y)) \frac{\partial \varphi}{\partial x} \de x \geq \int f \Phi(1,0) = \Phi\left(\int f,\int g\right). $$
A similar argument shows that the result holds if $\int f = + \infty$ or $\int g = + \infty$. Thus we assume thereafter that $0 < \int f < +\infty$ and $0 < \int g < +\infty$.

Let us show that one may assume that $\int f = \int g = 1$. Let us define, for $x,y \in \R$ and $a,b \geq 0$,
$$ \widetilde{f}(x) = f\left(\Phi\left(\int f, 0 \right) x \right) \Phi(1,0), \quad \widetilde{g}(x) = g\left(\Phi \left(0, \int g \right) x \right) \Phi(0,1), $$
$$ \widetilde{h}(x) = h\left(\Phi \left(\int f, \int g \right) x \right), $$
$$ \widetilde{\varphi}(x,y) = \frac{\varphi(\Phi(\int f, 0) x, \Phi(0, \int g) y)}{\Phi(\int f, \int g)}, \quad \widetilde{\Phi}(a,b) = \Phi\left( a \frac{\int f}{\Phi(\int f, \int g)}, b \frac{\int g}{\Phi(\int f, \int g)} \right). $$
Let $x \in \supp(\widetilde{f})$, $y \in \supp(\widetilde{g})$, and let $\widetilde{\rho}, \widetilde{\eta} > 0$. One has,
\begin{eqnarray*}
\widetilde{h} (\widetilde{\varphi}(x,y)) \left( \frac{\partial \widetilde{\varphi}}{\partial x} \widetilde{\rho} + \frac{\partial \widetilde{\varphi}}{\partial y} \widetilde{\eta} \right) & \geq & \Phi\left( f(\Phi(\int f, 0) x) \frac{\Phi(\int f, 0)}{\Phi(\int f, \int g)} \widetilde{\rho}, g(\Phi(0, \int g) y) \frac{\Phi(0, \int g)}{\Phi(\int f, \int g)} \widetilde{\eta} \right) \\ & = & \widetilde{\Phi}(\widetilde{f}(x) \widetilde{\rho}, \widetilde{g}(y) \widetilde{\eta}).
\end{eqnarray*}
Notice that the functions $\widetilde{\varphi}$ and $\widetilde{\Phi}$ satisfy the same assumptions as the functions $\varphi$ and $\Phi$ respectively, and that $\int \widetilde{f} = \int \widetilde{g} = 1$. If the result holds for functions of integral one, then
$$ \int \widetilde{h}(w) \de w \geq \widetilde{\Phi}(1,1) = 1. $$
The change of variable $w = z/\Phi(\int f, \int g)$ leads us to
$$ \int h(z) \de z \geq \Phi \left(\int f, \int g \right). $$

Assume now that $\int f = \int g = 1$. By standard approximation, one may assume that $f$ and $g$ are compactly supported positive Lipschitz functions (relying on the fact that $\Phi$ is continuous and increasing in each coordinate, compare with \cite[page 343]{Barthe1}). Thus there exists a non-decreasing map $T : \supp(f) \to \supp(g)$ such that for every $x \in \supp(f)$,
$$ f(x) = g(T(x))T'(x), $$
see e.g. \cite{Barthe2}, \cite{V}. Since $T$ is non-decreasing and $\partial \varphi / \partial x, \partial \varphi / \partial y > 0$, the function $\Theta : \supp(f) \to \varphi(\supp(f), T(\supp(f)))$ defined by $\Theta(x) = \varphi(x, T(x))$ is bijective. Hence the change of variable $z = \Theta(x)$ is admissible and one has,
\begin{eqnarray*}
\int h(z) \de z \geq \int_{\supp(f)} h(\varphi(x,T(x))) \left( \frac{\partial \varphi}{\partial x} + \frac{\partial \varphi}{\partial y} T'(x) \right) \de x & \geq & \int_{\supp(f)} \Phi(f(x), g(T(x))T'(x)) \de x \\ & = & \int \Phi(f(x), f(x)) \de x.
\end{eqnarray*}
Using homogeneity of $\Phi$, one deduces that
$$ \int h \geq \Phi(1, 1) \int f(x) \de x = \Phi \left(\int f, \int g \right). $$

{\bf Step 2 :} (Tensorization) \\
Let $n$ be a positive integer and assume that Theorem \ref{Borell} holds in $\R^n$. Let $f,g,h,\varphi, \Phi$ satisfying the assumptions of Theorem \ref{Borell} in $\R^{n+1}$. Recall that the inequality
\begin{eqnarray}\label{assumption}
h(\varphi(x,y)) \Pi_{k=1}^{n+1} \left( \frac{\partial \varphi_k}{\partial x_k} \rho_k + \frac{\partial \varphi_k}{\partial y_k} \eta_k \right) \geq \Phi(f(x)\Pi_{k=1}^{n+1} \rho_k, g(y) \Pi_{k=1}^{n+1} \eta_k),
\end{eqnarray}
holds for every $x \in \supp(f), y \in \supp(g)$, and for every $\rho_1, \dots, \rho_{n+1}, \eta_1, \dots, \eta_{n+1} > 0$. Let us define, for $x_{n+1}, y_{n+1}, z_{n+1} \in \R$,
$$ F(x_{n+1}) = \int_{\R^n} f(x, x_{n+1}) \de x, \quad G(y_{n+1}) = \int_{\R^n} g(x, g_{n+1}) \de x, \quad H(z_{n+1}) = \int_{\R^n} h(x, z_{n+1}) \de x. $$
Since $\int f > 0, \int g > 0$, the support of $F$ and the support of $G$ are nonempty. Let $x_{n+1} \in \supp(F), y_{n+1} \in \supp(G)$, and let $\rho_{n+1}, \eta_{n+1} > 0$. Let us define, for $x, y, z \in \R^n$,
$$ f_{x_{n+1}}(x) = f(x, x_{n+1}) \rho_{n+1}, \quad g_{y_{n+1}}(y) = g(y, y_{n+1}) \eta_{n+1}, \quad \overline{\varphi}(x,y) = (\varphi_1(x_1,y_1), \dots, \varphi_n(x_n,y_n)), $$
$$ h_{\varphi_{n+1}}(z) = h(z, \varphi_{n+1}(x_{n+1}, y_{n+1}))\left( \frac{\partial \varphi_{n+1}}{\partial x_{n+1}} \rho_{n+1} + \frac{\partial \varphi_{n+1}}{\partial y_{n+1}} \eta_{n+1} \right). $$
Let $x \in \supp(f_{x_{n+1}}), y \in \supp(g_{y_{n+1}})$, and let $\rho_1, \dots, \rho_n, \eta_1, \dots, \eta_n > 0$. One has
\begin{eqnarray*}
h_{\varphi_{n+1}}(\overline{\varphi}(x,y)) \Pi_{k=1}^n \left( \frac{\partial \overline{\varphi_k}}{\partial x_k} \rho_k + \frac{\partial \overline{\varphi_k}}{\partial y_k} \eta_k \right) & = & h(\varphi(x,x_{n+1},y,y_{n+1})) \Pi_{k=1}^{n+1} \left( \frac{\partial \varphi_k}{\partial x_k} \rho_k + \frac{\partial \varphi_k}{\partial y_k} \eta_k \right) \\ & \geq & \Phi(f(x,x_{n+1})\Pi_{k=1}^{n+1} \rho_k, g(y,y_{n+1}) \Pi_{k=1}^{n+1} \eta_k) \\ & = & \Phi(f_{x_{n+1}}(x)\Pi_{k=1}^{n} \rho_k, g_{y_{n+1}}(y) \Pi_{k=1}^{n} \eta_k),
\end{eqnarray*}
where the inequality follows from inequality (\ref{assumption}). Hence, applying Theorem \ref{Borell} in dimension $n$, one has
$$ \int_{\R^n} h_{\varphi_{n+1}}(x) \de x \geq \Phi\left( \int_{\R^n} f_{x_{n+1}}(x) \de x, \int_{\R^n} g_{y_{n+1}}(x) \de x \right). $$
This yields that for every $x_{n+1} \in \supp(F), y_{n+1} \in \supp(G)$, and for every $\rho_{n+1}, \eta_{n+1} > 0$,
$$ H(\varphi_{n+1}(x_{n+1}, y_{n+1})) \left( \frac{\partial \varphi_{n+1}}{\partial x_{n+1}} \rho_{n+1} + \frac{\partial \varphi_{n+1}}{\partial y_{n+1}} \eta_{n+1} \right) \geq \Phi(F(x_{n+1}), G(y_{n+1})). $$
Hence, applying Theorem \ref{Borell} in dimension 1, one has
$$ \int_{\R} H(x) \de x \geq \Phi\left( \int_{\R} F(x) \de x, \int_{\R} G(x) \de x \right). $$
This yields the desired inequality.
\end{proof}

\section{Applications to the log-Brunn-Minkowski inequality}

In this section, we discuss applications of the above inequalities to the log-Brunn-Minkowski inequality of B\"or\"oczky, Lutwak, Yang and Zhang \cite{BLYZ}. 

Recall that a {\it convex body} in $\R^n$ is a compact convex subset of $\R^n$ with nonempty interior. B\"or\"oczky et al. conjectured the following inequality.

\begin{conjecture}[log-Brunn-Minkowski inequality]\label{log-BM}

Let $K,L$ be symmetric convex bodies in $\R^n$ and let $\lambda \in [0,1]$. Then,
$$ |(1-\lambda) \cdot K \oplus_0 \lambda \cdot L| \geq |K|^{1-\lambda} |L|^{\lambda}. $$

\end{conjecture}

Here,
$$ (1-\lambda) \cdot K \oplus_0 \lambda \cdot L = \{ x \in \R^n : \langle x,u \rangle \leq h_K(u)^{1-\lambda} h_L(u)^{\lambda}, \mathrm{~for~all~} u\in S^{n-1} \}, $$
where $S^{n-1}$ denotes the $n$-dimensional Euclidean unit sphere, $h_K$ denotes the support function of $K$, defined by $h_K(u) = \max_{x \in K} \langle x,u \rangle$, and $|\cdot|$ stands for Lebesgue measure. 

B\"or\"oczky et al. \cite{BLYZ} proved that Conjecture \ref{log-BM} holds in the plane. Using Corollary \ref{B2} with ${\bf p}=(0, \dots, 0)$, Saroglou \cite{S1} proved that Conjecture \ref{log-BM} holds for unconditional convex bodies in $\R^n$ (a set $K \subset \R^n$ is {\it unconditional} if for every $(x_1, \dots, x_n) \in K$ and for every $(\eps_1, \dots, \eps_n) \in \{-1,1\}^n$, one has $(\eps_1 x_1, \dots, \eps_n x_n ) \in K$).

Recall that a measure $\mu$ is {\it $s$-concave}, $s \in [-\infty, +\infty]$, if the inequality
$$ \mu((1-\lambda) A + \lambda B) \geq M_s^{\lambda} (\mu(A), \mu(B)) $$
holds for all compact sets $A,B \subset \R^n$ such that $\mu(A)\mu(B) > 0$ and for every $\lambda \in [0,1]$ (see \cite{B1}, \cite{B2}). The $0$-concave measures are also called {\it log-concave measures}, and the $-\infty$-concave measures are also called {\it convex measures}. A function $f : \R^n \to [0, +\infty)$ is $\alpha$-concave, $\alpha \in [-\infty, +\infty]$, if the inequality
$$ f((1-\lambda) x + \lambda y) \geq M_{\alpha}^{\lambda} (f(x), f(y)) $$
holds for every $x,y \in \R^n$ such that $f(x)f(y) > 0$ and for every $\lambda \in [0,1]$.

Saroglou \cite{S2} recently proved that if the log-Brunn-Minkowski inequality holds, then the inequality
$$ \mu((1-\lambda) \cdot K \oplus_0 \lambda \cdot L) \geq \mu(K)^{1-\lambda} \mu(L)^{\lambda} $$
holds for every symmetric log-concave measure $\mu$, for all symmetric convex bodies $K,L$ in $\R^n$ and for every $\lambda \in [0,1]$.

An extension of the log-Brunn-Minkowski inequality for convex measures was proposed by the author in \cite{M4}, and reads as follows:

\begin{conjecture}\label{Lp-BM}

Let $p \in [0,1]$. Let $\mu$ be a symmetric measure in $\R^n$ that has an $\alpha$-concave density function, with $\alpha \geq - \frac{p}{n}$. Then for every symmetric convex body $K,L$ in $\R^n$ and for every $\lambda \in [0,1]$,
\begin{eqnarray}
\mu((1-\lambda) \cdot K \oplus_p \lambda \cdot L) \geq M_{\left(\frac{n}{p} + \frac{1}{\alpha}\right)^{-1}}^{\lambda} (\mu(K), \mu(L)).
\end{eqnarray}

\end{conjecture}

Here,
$$ (1-\lambda) \cdot K \oplus_p \lambda \cdot L = \{ x \in \R^n : \langle x,u \rangle \leq M_p^{\lambda}(h_K(u),h_L(u)), \mathrm{~for~all~} u\in S^{n-1} \}. $$

In Conjecture \ref{Lp-BM}, if $\alpha$ or $p$ is equal to $0$, then $(n/p + 1/\alpha)^{-1}$ is defined by continuity and is equal to $0$. Notice that Conjecture \ref{log-BM} is a particular case of Conjecture \ref{Lp-BM} when taking $\mu$ to be Lebesgue measure and $p=0$.

By using Corollary \ref{line}, we will prove that Conjecture \ref{log-BM} implies Conjecture \ref{Lp-BM}, when $\alpha \leq 1$, generalizing Saroglou's result discussed earlier.

\begin{theorem}\label{equiv}

If the log-Brunn-Minkowski inequality holds, then the inequality
$$ \mu((1-\lambda) \cdot K \oplus_p \lambda \cdot L) \geq M_{\left(\frac{n}{p} + \frac{1}{\alpha}\right)^{-1}}^{\lambda} (\mu(K), \mu(L)) $$
holds for every $p \in [0,1]$, for every symmetric measure $\mu$ in $\R^n$ that has an $\alpha$-concave density function, with $1 \geq \alpha \geq - \frac{p}{n}$, for every symmetric convex body $K,L$ in $\R^n$ and for every $\lambda \in [0,1]$.

\end{theorem}

\begin{proof}
Let $K_0, K_1$ be symmetric convex bodies in $\R^n$ and let $\lambda \in (0,1)$. Let us denote $K_{\lambda} = (1-\lambda) \cdot K_0 \oplus_p \lambda \cdot K_1$ and let us denote by $\psi$ the density function of $\mu$. Let us define, for $t > 0$, $h(t) = |K_{\lambda} \cap \{ \psi \geq t \}|$, $f(t) = |K_0 \cap \{ \psi \geq t \}|$ and $g(t) = |K_1 \cap \{ \psi \geq t \}|$. Notice that
$$ \mu(K_{\lambda}) = \int_{K_{\lambda}} \psi(x) \de x = \int_{K_{\lambda}} \int_0^{\psi(x)} \de t \de x = \int_0^{+\infty} |K_{\lambda} \cap \{ \psi \geq t \}| = \int_0^{+\infty} h(t) \de t. $$
Similarly, one has
$$ \mu(K_0) = \int_0^{+\infty} f(t) \de t, \qquad \mu(K_1) = \int_0^{+\infty} g(t) \de t. $$
Let $t,s > 0$ such that the sets $\{\psi \geq t\}$ and $\{\psi \geq s\}$ are nonempty. Let us denote $L_0 = \{\psi \geq t\}$, $L_1 = \{\psi \geq s\}$ and $L_{\lambda} = \{\psi \geq M_{\alpha}^{\lambda}(t,s) \}$. If $x \in L_0$ and $y \in L_1$, then $\psi((1-\lambda)x + \lambda y) \geq M_{\alpha}^{\lambda}(\psi(x),\psi(y)) \geq M_{\alpha}^{\lambda}(t,s)$. Hence,
$$ L_{\lambda} \supset (1-\lambda) L_0 + \lambda L_1 \supset (1-\lambda) \cdot L_0 \oplus_p \lambda \cdot L_1, $$
the last inclusion following from the fact that $p \leq 1$. We deduce that
$$ K_{\lambda} \cap L_{\lambda} \supset ((1-\lambda) \cdot K_0 \oplus_p \lambda \cdot K_1) \cap ((1-\lambda) \cdot L_0 \oplus_p \lambda \cdot L_1) \supset (1-\lambda) \cdot (K_0 \cap L_0) \oplus_p \lambda \cdot (K_1 \cap L_1). $$
Hence,
$$ h(M_{\alpha}^{\lambda}(t,s)) = |K_{\lambda} \cap L_{\lambda}| \geq |(1-\lambda) \cdot (K_0 \cap L_0) \oplus_p \lambda \cdot (K_1 \cap L_1)| \geq M_{\frac{p}{n}}^{\lambda}(f(t),g(s)), $$
the last inequality is valid for $p \geq 0$ and follows from the log-Brunn-Minkowski inequality by using homogeneity of Lebesgue measure (see \cite[beginning of section 3]{BLYZ}). Thus we may apply Corollary \ref{line} to conclude that
$$ \mu(K_{\lambda}) = \int_0^{+\infty} h \geq M_{\left(\frac{n}{p} + \frac{1}{\alpha}\right)^{-1}}^{\lambda}\left( \int_0^{+\infty} f, \int_0^{+\infty} g \right) = M_{\left(\frac{n}{p} + \frac{1}{\alpha}\right)^{-1}}^{\lambda}( \mu(K_0), \mu(K_1)). $$
\end{proof}

Since the log-Brunn-Minkowski inequality holds true in the plane, we deduce that Conjecture \ref{Lp-BM} holds true in the plane (with the restriction $\alpha \leq 1$). Notice that Conjecture \ref{Lp-BM} holds true in the unconditional case as a consequence of Corollary \ref{B2} (see \cite{M4}).

\vspace{1cm}

\noindent Arnaud Marsiglietti \\
Institute for Mathematics and its Applications \\
University of Minnesota \\
207 Church Street SE, 434 Lind Hall, \\
Minneapolis, MN 55455, USA \\
E-mail address: arnaud.marsiglietti@ima.umn.edu

\end{document}